# Martingale Representation and Logarithmic-Sobolev Inequality for Fractional Ornstein-Uhlenbeck Measure

Xiaoxia Sun[a,*], Feng Guo[b]

[a] School of Mathematics,
Dongbei University of Finance and Economics, Dalian, 116025, China
[b] School of Mathematical Sciences,
Dalian University of Technology, Dalian, 116024, China

## Abstract

In this paper, we consider the measure determined by a fractional Ornstein-Uhlenbeck process. For such measure, we establish a martingale representation theorem and consequently obtain the Logarithmic-Sobolev inequality. To this end, we also present the integration by parts formula for such measure, which is obtained via its pull back formula and the Bismut method.

## 1. Introduction

A stochastic diffusion process $(X_t)_{0 \le t \le 1}$ is called a fractional Ornstein-Uhlenbeck process if it satisfies the following stochastic differential equation

$$dX_t = -\alpha X_t dt + dB_t^H, \quad X_0 = 0. \tag{1}$$

where $\alpha > 0$ and $(B_t^H)_{0 \le t \le 1}$ is a $n$-dimensional fractional Brownian motion with Hurst parameter $H > 1/2$. Fractional Ornstein-Uhlenbeck processes are widely used to describe the long memory property for some time series, such as finance, hydrology, telecommunications, insurance and computer networks. The measure determined by a fractional Ornstein-Uhlenbeck process is called fractional Ornstein-Uhlenbeck measure. In this paper, we investigate the martingale representation and the Logarithmic-Sobolev inequality for such measure.

Quite a lot of interest has been attracted to the study of the martingale representation and Logarithmic-Sobolev inequalities for different measures. It is proved in [13] that Logarithmic-Sobolev inequality holds for Wiener measure on the path space over connected Lie group. For the Wiener measure on the path space over Riemannian manifold, [17] gives the Logarithmic-Sobolev inequality with a bounding constant which can be estimated in terms of the Ricci curvature. Moreover, for such measure, the Logarithmic-Sobolev inequality can also be obtained by embedding the manifold into an Euclidean space [3] and the martingale representation [5]. For Brwonian bridge measure on loop space, the martingale representation and the Logarithmic-Sobolev inequality are investigated in [2, 12, 14].

The Logarithmic-Sobolev inequalities for measures can be obtained via their martingale representations which can be established by their integration by parts formulas. The integration by parts formulas for different measures are important in infinite dimension analysis and have been well studied in the literature. For instance, the integration by parts formula is investigated for Wiener measure on the path space in [7, 11, 15], for Brownian bridge measure on the loop space in [1, 8, 10, 16], for fractional Wiener measure under different integrals in [6, 9] and for fractional Ornstein-Uhlenbeck measure in [21].

---
*Corresponding author
*Email addresses:* fly428@163.com (Xiaoxia Sun), fguo@dlut.edu.cn (Feng Guo)

**Contributions.** We establish the pull back formula (Proposition 3.1) and an integration by parts formula for fractional Ornstein-Uhlenbeck measure (Theorem 3.2). We give a martingale representation theorem (Theorem 4.1) by the corresponding integration by parts formula. Consequently, we derive the Logarithmic-Sobolev inequality (Theorem 4.2) by the martingale representation theorem.

The paper is organized as follows. In Section 2, we give some preliminaries about fractional Brownian motions. We present in Section 3 the pull back formula and the integration by parts formula. In Section 4, we obtain the martingale representation theorem and the Logarithmic-Sobolev inequality for fractional Ornstein-Uhlenbeck measure.

## 2. Preliminaries

Let $(\Omega, \mathcal{F}, \mathcal{F}_t, \nu)$ be a filtered probability space, where $\Omega$ is the space of $\mathbb{R}^n$-valued continuous functions on $[0, 1]$ with the initial value zero and the topology of local uniform convergence, $\nu$ is the fractional Ornstein-Uhlenbeck measure such that coordinate process $(X_t(\omega))_{0 \leq t \leq 1} = (\omega_t)_{0 \leq t \leq 1}$ satisfies (1), $\mathcal{F}$ is the $\nu$-completion of the Borel $\sigma$-algebra of $\Omega$ and $\mathcal{F}_t$ is the $\nu$-completed natural filtration of $\omega$. In the following we consider fractional integral with $H > 1/2$.

For any $p \in [1, \infty)$, let

$$L^p(\Omega; \nu) = \{F \mid F : \Omega \to \mathbb{R}, \|F\|_p := (\mathbb{E}_\nu |F|^p)^{\frac{1}{p}} < \infty\}.$$

We denote $I_{0+}^{H+1/2}(L^2(\Omega; \nu))$ as $(H + 1/2)$-Hölder left fractional Riemann-Liouville integral operator. By [6], the isomorphism operator $K : L^2(\Omega; \nu) \to I_{0+}^{H+1/2}(L^2(\Omega; \nu))$ is defined by

$$(Kh)_t = \int_0^t K(t, s) h_s ds,$$

where $h \in L^2(\Omega; \nu)$ and

$$K(t, s) = c_H s^{\frac{1}{2} - H} \int_s^t u^{H - \frac{1}{2}} (u - s)^{H - \frac{3}{2}} du I_{[0, t]}(s), \tag{2}$$

in which for beta function $B(\cdot, \cdot)$,

$$c_H = \sqrt{\frac{H(2H - 1)}{B(2 - 2H, H - \frac{1}{2})}}. \tag{3}$$

By the definition of $Kh$, for $H > 1/2$, we have

$$\begin{aligned}
(K^{-1} h)_t &= t^{H - \frac{1}{2}} D_{0+}^{H - \frac{1}{2}} t^{\frac{1}{2} - H} h' \\
&= \frac{1}{\Gamma(\frac{3}{2} - H)} \left( t^{\frac{1}{2} - H} h_t' + (H - \frac{1}{2}) t^{H - \frac{1}{2}} \int_0^t \frac{t^{\frac{1}{2} - H} h_t' - u^{\frac{1}{2} - H} h_u'}{(t - u)^{\frac{1}{2} + H}} du \right),
\end{aligned} \tag{4}$$

where $D_{0+}^{H-1/2}$ is $(H - 1/2)$-order left-sided Riemann-Liouville derivative and $h'$ is the derivative of $h$. By [6, Theorem 3.3], the Cameron-Martin vector field on $\Omega$ is

$$\mathcal{H}_H = \{Kh \mid h \text{ is adapted process and } h \in L^2(\Omega; \nu)\}.$$

The scalar product on $\mathcal{H}_H$ is defined by

$$\langle Kh, Kg \rangle_{\mathcal{H}_H} = \langle h, g \rangle_{L^2(\Omega; \nu)} = \mathbb{E}_\nu \left[ \int_0^1 \langle h_t, g_t \rangle dt \right].$$

For $Kh \in \mathcal{H}_H$, the directional derivative of $F$ along $Kh$ is

$$D_h F(\omega) = \lim_{\delta \to 0} \frac{1}{\delta} \left( F\left(\omega + \delta(Kh)\right) - F(\omega) \right).$$

We denote all the smooth cylindrical functions on $\Omega$ by

$$\mathcal{F}C^\infty(\Omega) = \{ F \mid F(\omega) = f(\omega_{t_1}, ..., \omega_{t_n}), 0 < t_1 \le t_2 \le \cdots \le t_n \le 1, f \in C^\infty(\mathbb{R}^n) \}.$$

For $F \in \mathcal{F}C^\infty(\Omega)$, the directional derivative of $F$ along $Kh$ is

$$D_h F(\omega) = \sum_{i=1}^n \langle \nabla^i F, Kh \rangle_{\mathbb{R}^n},$$

where

$$\nabla^i F = \nabla^i f(\omega_{t_1}, \cdots, \omega_{t_n})$$

is the gradient with respect to the $i$ variable of $f$. The gradient $DF : \Omega \to \mathcal{H}_H$ is determined by

$$\langle DF, Kh \rangle_{\mathcal{H}_H} = D_h F.$$

We denote the domain of $D$ by $Dom(D)$.

## 3. Integration by parts formula for fractional Ornstein-Uhlenbeck measure

To obtain the integration by parts formula for the fractional Ornstein-Uhlenbeck measure, inspired by the idea in [4], we first construct a $\mathbb{R}^n$-valued function $(\beta_t)_{0 \le t \le 1}$ such that for any $r \in (-\epsilon, \epsilon)$, the following stochastic differential equation

$$dX_t(r) = -\alpha X_t(r) dt + dB_t^H(r), \tag{5}$$

has a solution $(X_t(r))_{0 \le t \le 1}$ satisfying

Condition 1. $(X_t(r))_{0 \le t \le 1} \in \Omega$ for any $r$,

Condition 2. $\frac{d}{dr} X_t(r)\big|_{r=0}$ exists and $(Kh)_t = \frac{d}{dr} X_t(r)\big|_{r=0}$ for $(h_t)_{0 \le t \le 1} \in L^2(\Omega; \nu)$,

where $B_t^H(r)$ is defined by

$$B_t^H(r) = B_t^H + r \int_0^t \beta_s ds.$$

Note that $(X_t(0))_{0 \le t \le 1} = (X_t)_{0 \le t \le 1}$ and $(B_t^H(0))_{0 \le t \le 1} = (B_t^H)_{0 \le t \le 1}$.

**Proposition 3.1.** *If $(\beta_t)_{0 \le t \le 1}$ satisfies* Condition *1 and 2, then we have*

$$\beta_t = (Kh)_t' + \alpha(Kh)_t. \tag{6}$$

*Proof.* Differentiating (5) with respect to $r$ at $r = 0$, we get

$$d \frac{d}{dr} X_t(r)\big|_{r=0} = -\alpha \frac{d}{dr} X_t(r)\big|_{r=0} dt + d \frac{d}{dr} B_t^H(r)\big|_{r=0}.$$

By Condition 2, we have

$$d \frac{d}{dr} X_t(r)\big|_{r=0} = (Kh)_t' dt.$$

Then, it holds that

$$(Kh)_t' dt = -\alpha(Kh)_t dt + \beta_t dt,$$

which yields that $\beta_t = (Kh)_t' + \alpha(Kh)_t$. □

In the following we establish the integration by parts formula for the fractional Ornstein-Uhlenbeck measure $\nu$ via the pull back formula given in Proposition 3.1.

**Theorem 3.2.** *For $F \in Dom(D)$ and $Kh \in \mathcal{H}_H$, the integration by parts formula for fractional Ornstein-Uhlenbeck measure $\nu$ is*

$$\mathbb{E}_\nu \left[ F \int_0^1 \left\langle \left( K^{-1} \int_0^\cdot \beta_u du \right)_t, dB_t \right\rangle \right] = \mathbb{E}_\nu[D_h F], \tag{7}$$

*where $\beta_t = (Kh)_t' + \alpha(Kh)_t$.*

*Proof.* According to [6], there is a $n$-dimensional Brownian motion $(B_t)_{0 \leq t \leq 1}$ such that $B_t^H = \int_0^t K(t,s) dB_s$. Thus, by Proposition 3.1, we obtain

$$B_t^H(r) = \int_0^t K(t,s) d \left( B_s + r \int_0^s \left( K^{-1} \int_0^\cdot \beta_u du \right)_v dv \right).$$

For $t \in [0,1]$, we set

$$\rho_t = \exp \left\{ -r \int_0^t \left\langle \left( K^{-1} \int_0^\cdot \beta_u du \right)_s, dB_s \right\rangle - \frac{r^2}{2} \int_0^t \left( K^{-1} \int_0^\cdot \beta_u du \right)_s^2 ds \right\}. \tag{8}$$

For $H > 1/2$, by (4), we have

$$\begin{aligned}
&\left( K^{-1} \int_0^\cdot \beta_u du \right)_s \\
=& h_s + \alpha \left( K^{-1} \int_0^\cdot (Kh)_u du \right)_s \\
=& h_s + \frac{\alpha}{\Gamma(\frac{3}{2} - H)} \left( s^{\frac{1}{2} - H}(Kh)_s + (H - \frac{1}{2}) s^{H - \frac{1}{2}} \int_0^s \frac{s^{\frac{1}{2} - H}(Kh)_s - u^{\frac{1}{2} - H}(Kh)_u}{(s-u)^{\frac{1}{2} + H}} du \right) \\
=& h_s + \frac{\alpha s^{\frac{1}{2} - H}}{\Gamma(\frac{3}{2} - H)} (Kh)_s + \frac{\alpha(H - \frac{1}{2}) s^{H - \frac{1}{2}}}{\Gamma(\frac{3}{2} - H)} \int_0^s \frac{s^{\frac{1}{2} - H} - u^{\frac{1}{2} - H}}{(s-u)^{\frac{1}{2} + H}} (Kh)_u du \\
&+ \frac{\alpha(H - \frac{1}{2})}{\Gamma(\frac{3}{2} - H)} \int_0^s \frac{(Kh)_s - (Kh)_u}{(s-u)^{\frac{1}{2} + H}} du \\
=& h_s + I_1 + I_2 + I_3.
\end{aligned} \tag{9}$$

By Hölder inequality,

$$|(Kh)_s| \leq \left( \int_0^s K^2(s,u) du \right)^{\frac{1}{2}} \left( \int_0^1 |h|_u^2 du \right)^{\frac{1}{2}} \leq \frac{C_1 s^{1-H}}{\sqrt{2 - 2H}} \left( \int_0^1 |h|_u^2 du \right)^{\frac{1}{2}},$$

where $C_1$ satisfies

$$|K(s,t)| \leq C_1 t^{\frac{1}{2} - H}, \tag{10}$$

which is due to Theorem 3.2 in [6]. Hence, for $s \in [0,1]$, we have

$$|I_1| \leq \frac{\alpha C_1 s^{\frac{3}{2} - 2H}}{\Gamma(\frac{3}{2} - H)\sqrt{2 - 2H}} \left( \int_0^1 |h|_u^2 du \right)^{\frac{1}{2}}. \tag{11}$$

By [19], there is a constant $C_2$ such that

$$\int_0^s \frac{s^{\frac{1}{2} - H} - u^{\frac{1}{2} - H}}{(s-u)^{\frac{1}{2} + H}} du = C_2 s^{1 - 2H}, \tag{12}$$

then

$$|I_2| \leq \frac{-\left(H - \frac{1}{2}\right)\alpha C_1 C_2 s^{\frac{3}{2} - 2H}}{\Gamma(\frac{3}{2} - H)\sqrt{2 - 2H}} \left(\int_0^1 |h|_u^2 du\right)^{\frac{1}{2}}. \tag{13}$$

Since $Kh \in I_{0^+}^{H + 1/2}$, by [20, Theorem 3.6], $Kh$ is $H$-Hölder continuous. Therefore, there exists a constant $C_3$ such that

$$|I_3| \leq \frac{C_3 \alpha \left(H - \frac{1}{2}\right)}{\Gamma(\frac{3}{2} - H)} \int_0^s \frac{(s - u)^H}{(s - u)^{\frac{1}{2} + H}} du \left(\int_0^1 |h|_u^2 du\right)^{\frac{1}{2}} = \frac{C_3 \alpha \left(H - \frac{1}{2}\right)}{2\Gamma(\frac{3}{2} - H)} \left(\int_0^1 |h|_u^2 du\right)^{\frac{1}{2}}. \tag{14}$$

By (9), (11), (13) and (14), we conclude that

$$\begin{aligned}
&\int_0^1 \left|\left(K^{-1} \int_0^\cdot \beta_u du\right)_s\right|^2 ds \\
\leq &4 \int_0^1 |h_s|^2 ds + 4 \int_0^1 |I_1|^2 ds + 4 \int_0^1 |I_2|^2 ds + 4 \int_0^1 |I_3|^2 ds \\
\leq &4 \left(1 + \left(\frac{\alpha C_1}{\Gamma(\frac{3}{2} - H)\sqrt{(2 - 2H)(4 - 4H)}}\right)^2 + \left(\frac{\left(H - \frac{1}{2}\right)\alpha C_1 C_2}{\Gamma(\frac{3}{2} - H)\sqrt{(2 - 2H)(4 - 4H)}}\right)^2 \right. \\
&\left. + \left(\frac{C_3 \alpha \left(H - \frac{1}{2}\right)}{2\Gamma(\frac{3}{2} - H)}\right)^2\right) \left(\int_0^1 |h|_u^2 du\right).
\end{aligned} \tag{15}$$

Suppose that $h$ is bounded adapted process, then by (15), $\mathbb{E}_\nu[\rho_1] = 1$. It is easy to know that

$$\int_0^\cdot \beta_u du \in I_{0^+}^{H + \frac{1}{2}}(L^2(\Omega; \nu)).$$

Therefore, by [19, Theorem 2],

$$\int_0^t K(t, s) d\left(B_s + \left(K^{-1} \int_0^\cdot \beta_u du\right)_s\right)$$

is a fractional Brownian motion on $[0.1]$ under $\rho_1 \nu$. Hence, $(B_t^H(r))_{0 \leq t \leq 1}$ and $(B_t^H)_{0 \leq t \leq 1}$ have the same distribution under $\rho_1 \nu$ and $\rho_1$ respectively, which implies that $(X_t(r))_{0 \leq t \leq 1}$ and $(X_t)_{0 \leq t \leq 1}$ have the same distribution under $\rho_1 \nu$ and $\rho_1$ respectively. Therefore, for $F = f(X_{t_1}, ..., X_{t_n}) \in \mathcal{FC}^\infty(\Omega)$,

$$\mathbb{E}_{\rho_1 \nu}[f(X_{t_1}(r), \cdots, X_{t_n}(r))] = \mathbb{E}_\nu[f(X_{t_1}, \cdots, X_{t_n})].$$

Differentiating above equation with respect to $r$, we obtain

$$\begin{aligned}
&\frac{d}{dr} \mathbb{E}_\nu[\rho_1 f(X_{t_1}(r), \cdots, X_{t_n}(r))]\Big|_{r=0} \\
= &\mathbb{E}_\nu\left[\frac{d}{dr}\rho_1\Big|_{r=0} f(X_{t_1}, \cdots, X_{t_n})\right] + \mathbb{E}_\nu\left[\frac{d}{dr}f(X_{t_1}(r), \cdots, X_{t_n}(r))\Big|_{r=0}\right] \\
= &-\mathbb{E}_\nu\left[F \int_0^1 \left\langle \left(K^{-1} \int_0^\cdot \beta_u du\right)_t, dB_t\right\rangle\right] + \mathbb{E}_\nu[D_h F] = 0.
\end{aligned}$$

Hence, for bounded adapted process $h$, we have the following integration by parts formula

$$\mathbb{E}_\nu\left[F \int_0^1 \left\langle \left(K^{-1} \int_0^\cdot \beta_u du\right)_t, dB_t\right\rangle\right] = \mathbb{E}_\nu[D_h F]. \tag{16}$$

By (15), it is easy to know that $\left(K^{-1}\int_0^{\cdot}\beta_u du\right)_{0\le t\le 1} \in L^2(\Omega;\nu)$ for adapted process $h \in L^2(\Omega;\nu)$. Thus, integration by parts formula (16) holds for any adapted process $h \in L^2(\Omega;\nu)$. Moreover, since $D$ is a closable operator, integration by parts formula (16) holds for any $F \in Dom(D)$. $\qquad\square$

Now we have established the integration by parts formula for fractional Ornstein-Uhlenbeck measure $\nu$. In next section we will study the martingale representation and the Logarithmic-Sobolev inequality for $\nu$ via the integration by parts formula using the idea similar as in [12, 18].

## 4. Martingale representation theorem and Logarithmic-Sobolev inequality

To prove the Logarithmic-Sobolev inequality, we generalize the classical Clark-Ocone martingale representation theorem for Wiener measure to the fractional Ornstein-Uhlenbeck measure $\nu$. The classical Clark-Ocone martingale representation theorem states that every martingale adapted to the filtration $\mathcal{F}_t$ is a stochastic integral with respect to Brownian motion $B$. Suppose that $F \in L^2(\Omega;\nu)$, then there exists a $\mathcal{F}_t$-predictable process $\eta$ such that

$$F = \mathbb{E}_\nu[F] + \int_0^1 \langle \eta_t, dB_t\rangle.$$

The following theorem gives the explicit form of $\eta$ for the fractional Ornstein-Uhlenbeck measure.

**Theorem 4.1.** *Suppose that $F \in Dom(D)$, there exists a $\mathcal{F}_t$-predictable process $(\eta_t)_{0\le t\le 1}$ such that*

$$F = \mathbb{E}_\nu[F] + \int_0^1 \langle \eta_t, dB_t\rangle,$$

*where*

$$\eta_t = \mathbb{E}_\nu\left[(K^{-1}DF)_t + \int_t^1 K(s,t)\left(\int_s^1 \alpha^2 e^{-\alpha u} e^{\alpha s} P_u du - \alpha P_s\right)ds \,\bigg|\, \mathcal{F}_t\right], \tag{17}$$

*in which*

$$P_t = A_t(K^{-1}DF)_t + \frac{1}{\Gamma(\frac{3}{2}-H)}\int_t^1 (H-\frac{1}{2})u^{H-\frac{1}{2}}\frac{-t^{\frac{1}{2}-H}}{(u-t)^{\frac{1}{2}+H}}(K^{-1}DF)_u du,$$

$$A_t = \frac{1}{\Gamma(\frac{3}{2}-H)}\left(t^{\frac{1}{2}-H} + (H-\frac{1}{2})\int_0^t \frac{1}{(t-u)^{\frac{1}{2}+H}}du\right).$$

*Proof.* By the definition of $D_h F$, we have

$$\mathbb{E}_\nu[D_h F] = \mathbb{E}_\nu[\langle DF, Kh\rangle_{\mathcal{H}^H}] = \mathbb{E}_\nu\left[\int_0^1 \langle (K^{-1}DF)_t, h_t\rangle dt\right]. \tag{18}$$

At the same time, by the integration by parts formula (7), it holds that

$$\mathbb{E}_\nu[D_h F] = \mathbb{E}_\nu\left[\int_0^1 \langle \eta_t, dB_t\rangle \int_0^1 \left\langle \left(K^{-1}\int_0^{\cdot}\beta_u du\right)_t, dB_t\right\rangle\right] = \mathbb{E}_\nu\left[\int_0^1 \left\langle \eta_t, \left(K^{-1}\int_0^{\cdot}\beta_u du\right)_t\right\rangle dt\right]. \tag{19}$$

Let

$$j_t = \left(K^{-1}\int_0^{\cdot}\beta_u du\right)_t.$$

Then,

$$(Kh)_t + \alpha\int_0^t (Kh)_s ds = (Kj)_t,$$

which implies that

$$\int_0^t (Kh)_s ds = e^{-\alpha t}\int_0^t e^{\alpha s}(Kj)_s ds.$$

Thus

$$h_t = \left(K^{-1}\left(-\alpha e^{-\alpha \cdot}\int_0^{\cdot} e^{\alpha u}(Kj)_u du + (Kj)_{\cdot}\right)\right)_t. \tag{20}$$

Combining (18), (19) and (20), we have

$$\mathbb{E}_\nu\left[\int_0^1 \left\langle (K^{-1}DF)_t, \left(K^{-1}\left(-\alpha e^{-\alpha \cdot}\int_0^{\cdot} e^{\alpha u}(Kj)_u du + (Kj)_{\cdot}\right)\right)_t\right\rangle dt\right] = \mathbb{E}_\nu\left[\int_0^1 \langle \eta_t, j_t\rangle dt\right]. \tag{21}$$

By (4),

$$\begin{aligned}
&\left(K^{-1}\left(-\alpha e^{-\alpha \cdot}\int_0^{\cdot} e^{\alpha u}(Kj)_u du + (Kj)_{\cdot}\right)\right)_t\\
=&\frac{1}{\Gamma(\frac{3}{2}-H)}\left(t^{\frac{1}{2}-H} + (H-\frac{1}{2})\int_0^t \frac{1}{(t-u)^{\frac{1}{2}+H}}du\right)\left(\alpha^2 e^{-\alpha t}\int_0^t e^{\alpha u}(Kj)_u du - \alpha(Kj)_t\right)\\
&+\frac{(H-\frac{1}{2})t^{H-\frac{1}{2}}}{\Gamma(\frac{3}{2}-H)}\int_0^t \frac{-u^{\frac{1}{2}-H}\left(\alpha^2 e^{-\alpha u}\int_0^u e^{\alpha s}(Kj)_s ds - \alpha(Kj)_u\right)}{(t-u)^{\frac{1}{2}+H}}du + j_t\\
=&A_t\delta_t + \frac{(H-\frac{1}{2})t^{H-\frac{1}{2}}}{\Gamma(\frac{3}{2}-H)}\int_0^t \frac{-u^{\frac{1}{2}-H}\delta_u}{(t-u)^{\frac{1}{2}+H}}du + j_t,
\end{aligned}$$

where

$$\begin{aligned}
\delta_t =&\alpha^2 e^{-\alpha t}\int_0^t e^{\alpha u}(Kj)_u du - \alpha(Kj)_t,\\
A_t =&\frac{1}{\Gamma(\frac{3}{2}-H)}\left(t^{\frac{1}{2}-H} + (H-\frac{1}{2})\int_0^t \frac{1}{(t-u)^{\frac{1}{2}+H}}du\right).
\end{aligned} \tag{22}$$

Hence, the left side of (21) can be written as

$$\begin{aligned}
&\mathbb{E}_\nu\left[\int_0^1 \left\langle (K^{-1}DF)_t, \left(K^{-1}\left(-\alpha e^{-\alpha \cdot}\int_0^{\cdot} e^{\alpha u}(Kj)_u du + (Kj)_{\cdot}\right)\right)_t\right\rangle dt\right]\\
=&\mathbb{E}_\nu\left[\int_0^1 \langle (K^{-1}DF)_t, j_t\rangle dt\right] + \mathbb{E}_\nu\left[\int_0^1 \langle (K^{-1}DF)_t, A_t\delta_t\rangle dt\right]\\
&+\mathbb{E}_\nu\left[\int_0^1 \left\langle (K^{-1}DF)_t, \frac{(H-\frac{1}{2})t^{H-\frac{1}{2}}}{\Gamma(\frac{3}{2}-H)}\int_0^t \frac{-u^{\frac{1}{2}-H}\delta_u}{(t-u)^{\frac{1}{2}+H}}du\right\rangle dt\right].
\end{aligned} \tag{23}$$

By calculation, the third term of above equation is

$$\begin{aligned}
&\mathbb{E}_\nu\left[\int_0^1 \left\langle (K^{-1}DF)_t, \frac{(H-\frac{1}{2})t^{H-\frac{1}{2}}}{\Gamma(\frac{3}{2}-H)}\int_0^t \frac{-u^{\frac{1}{2}-H}\delta_u}{(t-u)^{\frac{1}{2}+H}}du\right\rangle dt\right]\\
=&\mathbb{E}_\nu\left[\int_0^1 \left\langle\int_u^1 \frac{(H-\frac{1}{2})t^{H-\frac{1}{2}}}{\Gamma(\frac{3}{2}-H)}\frac{-u^{\frac{1}{2}-H}}{(t-u)^{\frac{1}{2}+H}}(K^{-1}DF)_t dt, \delta_u\right\rangle du\right]\\
=&\mathbb{E}_\nu\left[\int_0^1 \left\langle\int_t^1 \frac{(H-\frac{1}{2})u^{H-\frac{1}{2}}}{\Gamma(\frac{3}{2}-H)}\frac{-t^{\frac{1}{2}-H}}{(u-t)^{\frac{1}{2}+H}}(K^{-1}DF)_u du, \delta_t\right\rangle dt\right].
\end{aligned}$$

Then (23) equals

$$\mathbb{E}_\nu\left[\int_0^1 \left\langle (K^{-1}DF)_t, j_t \right\rangle dt\right] + \mathbb{E}_\nu\left[\int_0^1 \left\langle A_t(K^{-1}DF)_t - \int_t^1 \frac{(H-\frac{1}{2})u^{H-\frac{1}{2}}t^{\frac{1}{2}-H}}{\Gamma(\frac{3}{2}-H)(u-t)^{\frac{1}{2}+H}}(K^{-1}DF)_u du, \delta_t \right\rangle dt\right]$$

$$=\mathbb{E}_\nu\left[\int_0^1 \left\langle (K^{-1}DF)_t, j_t \right\rangle dt\right] + \mathbb{E}_\nu\left[\int_0^1 \left\langle P_t, \delta_t \right\rangle dt\right],$$

(24)

where

$$P_t = A_t(K^{-1}DF)_t - \int_t^1 \frac{(H-\frac{1}{2})u^{H-\frac{1}{2}}t^{\frac{1}{2}-H}}{\Gamma(\frac{3}{2}-H)(u-t)^{\frac{1}{2}+H}}(K^{-1}DF)_u du.$$

(25)

By (22), we infer that

$$\mathbb{E}_\nu\left[\int_0^1 \langle P_t, \delta_t \rangle dt\right] = \mathbb{E}_\nu\left[\int_0^1 \left\langle P_t, \alpha^2 e^{-\alpha t}\int_0^t e^{\alpha u}(Kj)_u du - \alpha(Kj)_t \right\rangle dt\right]$$

$$=\mathbb{E}_\nu\left[\int_0^1 \left\langle P_t, \alpha^2 e^{-\alpha t}\int_0^t e^{\alpha u}(Kj)_u du \right\rangle dt\right] - \mathbb{E}_\nu\left[\int_0^1 \langle P_t, \alpha(Kj)_t \rangle dt\right]$$

$$=\mathbb{E}_\nu\left[\int_0^1 \left\langle \int_u^1 \alpha^2 e^{-\alpha t}e^{\alpha u}P_t dt, (Kj)_u \right\rangle du\right] - \mathbb{E}_\nu\left[\int_0^1 \langle P_t, \alpha(Kj)_t \rangle dt\right]$$

$$=\mathbb{E}_\nu\left[\int_0^1 \left\langle \int_t^1 \alpha^2 e^{-\alpha u}e^{\alpha t}P_u du, (Kj)_t \right\rangle dt\right] - \mathbb{E}_\nu\left[\int_0^1 \langle P_t, \alpha(Kj)_t \rangle dt\right]$$

$$=\mathbb{E}_\nu\left[\int_0^1 \left\langle \int_t^1 \alpha^2 e^{-\alpha u}e^{\alpha t}P_u du - \alpha P_t, \int_0^t K(t,s)j_s ds \right\rangle dt\right]$$

$$=\mathbb{E}_\nu\left[\int_0^1 \left\langle \int_s^1 K(t,s)\left(\int_t^1 \alpha^2 e^{-\alpha u}e^{\alpha t}P_u du - \alpha P_t\right) dt, j_s \right\rangle ds\right]$$

$$=\mathbb{E}_\nu\left[\int_0^1 \left\langle \int_t^1 K(s,t)\left(\int_s^1 \alpha^2 e^{-\alpha u}e^{\alpha s}P_u du - \alpha P_s\right) ds, j_t \right\rangle dt\right].$$

(26)

Combining (21), (24) and (26), we obtain

$$\mathbb{E}_\nu\left[\int_0^1 \left\langle (K^{-1}DF)_t + \int_t^1 K(s,t)\left(\int_s^1 \alpha^2 e^{-\alpha u}e^{\alpha s}P_u du - \alpha P_s\right) ds, j_t \right\rangle dt\right] = \mathbb{E}_\nu\left[\int_0^1 \langle \eta_t, j_t \rangle dt\right],$$

which yields that

$$\eta_t = \mathbb{E}_\nu\left[(K^{-1}DF)_t + \int_t^1 K(s,t)\left(\int_s^1 \alpha^2 e^{-\alpha u}e^{\alpha s}P_u du - \alpha P_s\right) ds \bigg| \mathcal{F}_t\right].$$

Here ends the proof. $\qquad\square$

In the following we establish the Logarithmic-Sobolev inequality for $\nu$ by its martingale representation theorem (Theorem 4.1).

**Theorem 4.2.** *For $F \in Dom(D)$, it holds that*

$$\mathbb{E}_\nu[F^2 \ln F^2] \le 4\left(1 + \frac{4\alpha^4 e^{2\alpha}C_1^2(1+\alpha^2)C}{2-2H} + \frac{2\alpha^2\widehat{C}}{2-2H}\right)\langle DF, DF \rangle_{\mathcal{H}_H} + \mathbb{E}_\nu[F^2]\ln\mathbb{E}_\nu[F^2],$$

where

$$C = \frac{1}{\Gamma^2(\frac{3}{2}-H)(2-2H)} + \frac{4}{\Gamma^2(\frac{3}{2}-H)(2-2H)} + \frac{C_2^2(H-\frac{1}{2})^2}{\Gamma^2(\frac{3}{2}-H)(2-2H)},$$

$$\widehat{C} = \frac{C_1^2}{\Gamma^2(\frac{3}{2}-H)(2-2H)} + \frac{4C_1^2}{\Gamma^2(\frac{3}{2}-H)(2-2H)}$$

$$+ \frac{2c_H^2\left(\left(B(H-\frac{1}{2},\frac{3}{2}-H)\right)^2 + \frac{1}{(H-\frac{1}{2})^2}\right)}{\Gamma^2(\frac{3}{2}-H)} + \frac{C_1^2 C_2^2(H-\frac{1}{2})^2}{\Gamma^2(\frac{3}{2}-H)(2-2H)},$$

$C_1$ and $C_2$ satisfy (10) and (12), respectively.

*Proof.* Let $G = F^2$ and $G_t$ be a right continuous version of $\mathbb{E}_\nu[G|\mathcal{F}_t]$. Then by Theorem 4.1,

$$dG_t = \langle \eta_t, dB_t \rangle,$$

where $\eta$ satisfies (17) in which $F$ is replaced by $G$. By Itô formula, we obtain

$$d(G_t \ln(G_t)) = (1 + \ln(G_t))dG_t + \frac{1}{2}\frac{|\eta_t|^2}{G_t}dt$$

$$= \langle (1 + \ln(G_t))\eta_t, dB_t \rangle + \frac{1}{2}\frac{|\eta_t|^2}{G_t}dt.$$

Due to $G_1 = \mathbb{E}_\nu[G|\mathcal{F}_1] = G$ and $G_0 = \mathbb{E}_\nu[G|\mathcal{F}_0] = \mathbb{E}_\nu[G]$, we get

$$\mathbb{E}_\nu[G \ln G] - \mathbb{E}_\nu[G]\ln \mathbb{E}_\nu[G] = \frac{1}{2}\mathbb{E}_\nu\left[\int_0^1 \frac{|\eta_t|^2}{G_t}dt\right]. \tag{27}$$

Since $DF^2 = 2FDF$, it holds that

$$\eta_t = \mathbb{E}_\nu\left[2F\left((K^{-1}DF)_t + \int_t^1\left(K(s,t)\int_s^1 \alpha^2 e^{-\alpha u}e^{\alpha s}P_u du - \alpha P_s\right)ds\right)\bigg|\mathcal{F}_t\right].$$

It follows that

$$|\eta_t|^2 \leq 4\mathbb{E}_\nu\left[F^2|\mathcal{F}_t\right]\mathbb{E}_\nu\left[\left|(K^{-1}DF)_t + \int_t^1 K(s,t)\left(\int_s^1 \alpha^2 e^{-\alpha u}e^{\alpha s}P_u du - \alpha P_s\right)ds\right|^2\bigg|\mathcal{F}_t\right]$$

$$\leq 8\mathbb{E}_\nu\left[F^2|\mathcal{F}_t\right]\mathbb{E}_\nu\left[|(K^{-1}DF)_t|^2 + 2\left(\int_t^1 \alpha^2 K(s,t)e^{\alpha s}\int_s^1 e^{-\alpha u}P_u du ds\right)^2\right.$$

$$\left. + 2\alpha^2\left(\int_t^1 K(s,t)P_s ds\right)^2\bigg|\mathcal{F}_t\right] \tag{28}$$

$$\leq 8\mathbb{E}_\nu\left[F^2|\mathcal{F}_t\right]\mathbb{E}_\nu\left[|(K^{-1}DF)_t|^2 + 2\alpha^4 e^{2\alpha}\int_t^1 (K(s,t))^2 ds \int_t^1\left(\int_s^1 e^{-\alpha u}P_u du\right)^2 ds\right.$$

$$\left. + 2\alpha^2\left(\int_t^1 K(s,t)P_s ds\right)^2\bigg|\mathcal{F}_t\right],$$

which is due to Hölder inequality and the inequality $(a+b)^2 \leq 2a^2 + 2b^2$. It is obvious that

$$
\begin{aligned}
\int_t^1 \left( \int_s^1 e^{-\alpha u} P_u du \right)^2 ds &= \int_t^1 \left( -\int_s^1 e^{-\alpha u} d \int_u^1 P_v dv \right)^2 ds \\
&\leq \int_0^1 \left( e^{-\alpha s} \int_s^1 P_v dv - \alpha \int_s^1 \int_u^1 P_v dv e^{-\alpha u} du \right)^2 ds \\
&\leq 2 \int_0^1 \left( \int_s^1 P_v dv \right)^2 + \alpha^2 \left( \int_s^1 P_v dv \right)^2 ds.
\end{aligned}
\tag{29}
$$

Hence, by (10), (28) and (29), we gain

$$
\begin{aligned}
|\eta_t|^2 \leq &8 \mathbb{E}_\nu \left[ F^2 | \mathcal{F}_t \right] \mathbb{E}_\nu \left[ |(K^{-1}DF)_t|^2 + 4\alpha^4 e^{2\alpha} C_1^2 t^{1-2H} (1+\alpha^2) \int_0^1 \left( \int_s^1 P_v dv \right)^2 ds \right. \\
&\left. + 2\alpha^2 \left( \int_t^1 K(s,t) P_s ds \right)^2 \,\middle|\, \mathcal{F}_t \right].
\end{aligned}
\tag{30}
$$

In the following we estimate $\left( \int_t^1 P_s ds \right)^2$ and $\left( \int_t^1 K(s,t) P_s ds \right)^2$ for (30) in two steps.

(I) The expression of $P$ implies that

$$
\begin{aligned}
\int_t^1 P_s ds = &\int_t^1 \frac{1}{\Gamma(\frac{3}{2}-H)} s^{\frac{1}{2}-H} (K^{-1}DF)_s ds + \frac{(H-\frac{1}{2})}{\Gamma(\frac{3}{2}-H)} \int_t^1 \int_0^s \frac{1}{(s-u)^{\frac{1}{2}+H}} du (K^{-1}DF)_s ds \\
&+ \frac{(H-\frac{1}{2})}{\Gamma(\frac{3}{2}-H)} \int_t^1 \int_s^1 u^{H-\frac{1}{2}} \frac{-s^{\frac{1}{2}-H}}{(u-s)^{\frac{1}{2}+H}} (K^{-1}DF)_u du ds.
\end{aligned}
\tag{31}
$$

It is obvious that

$$
\begin{aligned}
&\frac{(H-\frac{1}{2})}{\Gamma(\frac{3}{2}-H)} \int_t^1 \int_s^1 u^{H-\frac{1}{2}} \frac{-s^{\frac{1}{2}-H}}{(u-s)^{\frac{1}{2}+H}} (K^{-1}DF)_u du ds \\
=&\frac{(H-\frac{1}{2})}{\Gamma(\frac{3}{2}-H)} \int_t^1 \int_t^u u^{H-\frac{1}{2}} \frac{-s^{\frac{1}{2}-H}}{(u-s)^{\frac{1}{2}+H}} (K^{-1}DF)_u ds du \\
=&\frac{(H-\frac{1}{2})}{\Gamma(\frac{3}{2}-H)} \int_t^1 \int_t^s s^{H-\frac{1}{2}} \frac{-u^{\frac{1}{2}-H}}{(s-u)^{\frac{1}{2}+H}} du (K^{-1}DF)_s ds.
\end{aligned}
$$

Then (31) becomes

$$
\begin{aligned}
\int_t^1 P_s ds = &\int_t^1 \frac{1}{\Gamma(\frac{3}{2}-H)} s^{\frac{1}{2}-H} (K^{-1}DF)_s ds + \frac{(H-\frac{1}{2})}{\Gamma(\frac{3}{2}-H)} \int_t^1 \int_0^t \frac{1}{(s-u)^{\frac{1}{2}+H}} du (K^{-1}DF)_s ds \\
&+ \frac{(H-\frac{1}{2})}{\Gamma(\frac{3}{2}-H)} \int_t^1 \int_t^s \frac{1}{(s-u)^{\frac{1}{2}+H}} du (K^{-1}DF)_s ds \\
&+ \frac{(H-\frac{1}{2})}{\Gamma(\frac{3}{2}-H)} \int_t^1 \int_t^s s^{H-\frac{1}{2}} \frac{-u^{\frac{1}{2}-H}}{(s-u)^{\frac{1}{2}+H}} du (K^{-1}DF)_s ds \\
= &\int_t^1 \frac{1}{\Gamma(\frac{3}{2}-H)} s^{\frac{1}{2}-H} (K^{-1}DF)_s ds + \frac{(H-\frac{1}{2})}{\Gamma(\frac{3}{2}-H)} \int_t^1 \int_0^t \frac{1}{(s-u)^{\frac{1}{2}+H}} du (K^{-1}DF)_s ds \\
&+ \frac{(H-\frac{1}{2})}{\Gamma(\frac{3}{2}-H)} \int_t^1 \int_t^s \frac{s^{\frac{1}{2}-H} - u^{\frac{1}{2}-H}}{(s-u)^{\frac{1}{2}+H}} du s^{H-\frac{1}{2}} (K^{-1}DF)_s ds.
\end{aligned}
\tag{32}
$$

It holds that

$$\left(\int_t^1 P_s ds\right)^2$$

$$\leq \frac{1}{\Gamma^2(\frac{3}{2}-H)(2-2H)}\int_0^1 |(K^{-1}DF)_s|^2 ds + \frac{(H-\frac{1}{2})^2}{\Gamma^2(\frac{3}{2}-H)}\int_t^1\left(\int_0^t \frac{1}{(s-u)^{\frac{1}{2}+H}}du\right)^2 ds\int_0^1 |(K^{-1}DF)_s|^2 ds$$

$$+ \frac{(H-\frac{1}{2})^2}{\Gamma^2(\frac{3}{2}-H)}\int_0^1\left(\int_t^s \frac{s^{\frac{1}{2}-H}-u^{\frac{1}{2}-H}}{(s-u)^{\frac{1}{2}+H}}du\, s^{H-\frac{1}{2}}\right)^2 ds\int_0^1 |(K^{-1}DF)_s|^2 ds. \tag{33}$$

Since

$$\int_0^t \frac{1}{(s-u)^{\frac{1}{2}+H}}du = \frac{1}{\frac{1}{2}-H}(s^{\frac{1}{2}-H}-(s-t)^{\frac{1}{2}-H}),$$

we get

$$\frac{(H-\frac{1}{2})^2}{\Gamma^2(\frac{3}{2}-H)}\int_t^1\left(\int_0^t \frac{1}{(s-u)^{\frac{1}{2}+H}}du\right)^2 ds \leq \frac{(H-\frac{1}{2})^2}{\Gamma^2(\frac{3}{2}-H)}\int_t^1 \frac{2}{\left(\frac{1}{2}-H\right)^2}(s^{1-2H}+(s-t)^{1-2H})ds$$

$$= \frac{2}{\Gamma^2(\frac{3}{2}-H)(2-2H)}(1+(1-t)^{2-2H}-t^{2-2H}) \tag{34}$$

$$\leq \frac{4}{\Gamma^2(\frac{3}{2}-H)(2-2H)}.$$

By (12), we obtain

$$\frac{(H-\frac{1}{2})^2}{\Gamma^2(\frac{3}{2}-H)}\int_0^1\left(\int_t^s \frac{s^{\frac{1}{2}-H}-u^{\frac{1}{2}-H}}{(s-u)^{\frac{1}{2}+H}}du\, s^{H-\frac{1}{2}}\right)^2 ds$$

$$\leq \frac{(H-\frac{1}{2})^2}{\Gamma^2(\frac{3}{2}-H)}\int_0^1\left(\int_0^s \frac{u^{\frac{1}{2}-H}-s^{\frac{1}{2}-H}}{(s-u)^{\frac{1}{2}+H}}du\right)^2 s^{2H-1} ds \tag{35}$$

$$= \frac{(H-\frac{1}{2})^2}{\Gamma^2(\frac{3}{2}-H)}\int_0^1 C_2^2 s^{2-4H}s^{2H-1}ds \leq \frac{C_2^2(H-\frac{1}{2})^2}{\Gamma^2(\frac{3}{2}-H)(2-2H)}.$$

Combining (33), (34) and (35), it holds that

$$\left(\int_t^1 P_s ds\right)^2 \leq C\int_0^1 |(K^{-1}DF)_s|^2 ds, \tag{36}$$

where

$$C = \frac{1}{\Gamma^2(\frac{3}{2}-H)(2-2H)} + \frac{4}{\Gamma^2(\frac{3}{2}-H)(2-2H)} + \frac{C_2^2(H-\frac{1}{2})^2}{\Gamma^2(\frac{3}{2}-H)(2-2H)}.$$

(II) By (25), we have

$$\int_t^1 K(s,t)P_s ds$$

$$= \int_t^1 \frac{1}{\Gamma(\frac{3}{2}-H)}s^{\frac{1}{2}-H}K(s,t)(K^{-1}DF)_s ds + \frac{(H-\frac{1}{2})}{\Gamma(\frac{3}{2}-H)}\int_t^1\int_0^s \frac{1}{(s-u)^{\frac{1}{2}+H}}du\, K(s,t)(K^{-1}DF)_s ds \tag{37}$$

$$+ \frac{(H-\frac{1}{2})}{\Gamma(\frac{3}{2}-H)}\int_t^1\int_s^1 u^{H-\frac{1}{2}}\frac{-s^{\frac{1}{2}-H}}{(u-s)^{\frac{1}{2}+H}}(K^{-1}DF)_u du\, K(s,t)ds.$$

It is obvious that

$$
\frac{(H-\frac{1}{2})}{\Gamma(\frac{3}{2}-H)}\int_t^1\int_s^1 u^{H-\frac{1}{2}}\frac{-s^{\frac{1}{2}-H}}{(u-s)^{\frac{1}{2}+H}}(K^{-1}DF)_u du K(s,t)ds
$$
$$
=\frac{(H-\frac{1}{2})}{\Gamma(\frac{3}{2}-H)}\int_t^1\int_t^u u^{H-\frac{1}{2}}\frac{-s^{\frac{1}{2}-H}}{(u-s)^{\frac{1}{2}+H}}(K^{-1}DF)_u K(s,t)ds du
$$
$$
=\frac{(H-\frac{1}{2})}{\Gamma(\frac{3}{2}-H)}\int_t^1\int_t^s s^{H-\frac{1}{2}}\frac{-u^{\frac{1}{2}-H}}{(s-u)^{\frac{1}{2}+H}}K(u,t)du(K^{-1}DF)_s ds.
$$

Hence, by (37),

$$
\int_t^1 K(s,t)P_s ds
$$
$$
=\int_t^1\frac{1}{\Gamma(\frac{3}{2}-H)}s^{\frac{1}{2}-H}K(s,t)(K^{-1}DF)_s ds+\frac{(H-\frac{1}{2})}{\Gamma(\frac{3}{2}-H)}\int_t^1\int_0^t\frac{1}{(s-u)^{\frac{1}{2}+H}}du K(s,t)(K^{-1}DF)_s ds
$$
$$
+\frac{(H-\frac{1}{2})}{\Gamma(\frac{3}{2}-H)}\int_t^1\int_t^s\frac{s^{\frac{1}{2}-H}K(s,t)-u^{\frac{1}{2}-H}K(u,t)}{(s-u)^{\frac{1}{2}+H}}du s^{H-\frac{1}{2}}(K^{-1}DF)_s ds
$$
$$
=\int_t^1\frac{1}{\Gamma(\frac{3}{2}-H)}s^{\frac{1}{2}-H}K(s,t)(K^{-1}DF)_s ds+\frac{(H-\frac{1}{2})}{\Gamma(\frac{3}{2}-H)}\int_t^1\int_0^t\frac{1}{(s-u)^{\frac{1}{2}+H}}du K(s,t)(K^{-1}DF)_s ds \tag{38}
$$
$$
+\frac{(H-\frac{1}{2})}{\Gamma(\frac{3}{2}-H)}\int_t^1\int_t^s\frac{s^{\frac{1}{2}-H}(K(s,t)-K(u,t))}{(s-u)^{\frac{1}{2}+H}}du s^{H-\frac{1}{2}}(K^{-1}DF)_s ds
$$
$$
+\frac{(H-\frac{1}{2})}{\Gamma(\frac{3}{2}-H)}\int_t^1\int_t^s\frac{(s^{\frac{1}{2}-H}-u^{\frac{1}{2}-H})K(u,t)}{(s-u)^{\frac{1}{2}+H}}du s^{H-\frac{1}{2}}(K^{-1}DF)_s ds.
$$

Then by Hölder inequality, we get

$$
\left(\int_t^1 K(s,t)P_s ds\right)^2
$$
$$
\leq\frac{C_1^2 t^{1-2H}}{\Gamma^2(\frac{3}{2}-H)(2-2H)}\int_0^1|(K^{-1}DF)_s|^2 ds
$$
$$
+\frac{(H-\frac{1}{2})^2}{\Gamma^2(\frac{3}{2}-H)}\int_t^1\left(\int_0^t\frac{1}{(s-u)^{\frac{1}{2}+H}}du K(s,t)\right)^2 ds\int_0^1|(K^{-1}DF)_s|^2 ds
$$
$$
+\frac{(H-\frac{1}{2})^2}{\Gamma^2(\frac{3}{2}-H)}\int_t^1\left(\int_t^s\frac{s^{\frac{1}{2}-H}(K(s,t)-K(u,t))}{(s-u)^{\frac{1}{2}+H}}du s^{H-\frac{1}{2}}\right)^2 ds\int_0^1|(K^{-1}DF)_s|^2 ds \tag{39}
$$
$$
+\frac{(H-\frac{1}{2})^2}{\Gamma^2(\frac{3}{2}-H)}\int_t^1\left(\int_t^s\frac{(s^{\frac{1}{2}-H}-u^{\frac{1}{2}-H})K(u,t)}{(s-u)^{\frac{1}{2}+H}}du s^{H-\frac{1}{2}}\right)^2 ds\int_0^1|(K^{-1}DF)_s|^2 ds.
$$

Similar to (34), it holds that

$$
\int_t^1\left(\int_0^t\frac{1}{(s-u)^{\frac{1}{2}+H}}du K(s,t)\right)^2 ds\leq\frac{2C_1^2 t^{1-2H}}{(\frac{1}{2}-H)^2(2-2H)}(1+(1-t)^{2-2H}-t^{2-2H}) \tag{40}
$$
$$
\leq\frac{4C_1^2 t^{1-2H}}{(\frac{1}{2}-H)^2(2-2H)}.
$$

By the expression of $K(s,t)$ in (2), we have

$$\left(\int_t^s \frac{s^{\frac{1}{2}-H}(K(s,t)-K(u,t))}{(s-u)^{\frac{1}{2}+H}}du s^{H-\frac{1}{2}}\right)^2$$

$$=\left(c_H t^{\frac{1}{2}-H}\int_t^s\int_u^s\frac{v^{H-\frac{1}{2}}(v-t)^{H-\frac{3}{2}}}{(s-u)^{\frac{1}{2}+H}}dvdu\right)^2$$

$$=\left(c_H t^{\frac{1}{2}-H}\int_t^s v^{H-\frac{1}{2}}(v-t)^{H-\frac{3}{2}}\int_t^v\frac{1}{(s-u)^{\frac{1}{2}+H}}dudv\right)^2$$

$$=\left(-\frac{c_H t^{\frac{1}{2}-H}}{\frac{1}{2}-H}\int_t^s v^{H-\frac{1}{2}}(v-t)^{H-\frac{3}{2}}((s-v)^{\frac{1}{2}-H}-(s-t)^{\frac{1}{2}-H})dv\right)^2 \quad (41)$$

$$\leq\frac{2c_H^2 t^{1-2H}}{(\frac{1}{2}-H)^2}\left(\left(\int_t^s v^{H-\frac{1}{2}}(v-t)^{H-\frac{3}{2}}(s-v)^{\frac{1}{2}-H}dv\right)^2\right.$$

$$\left.+\left(\int_t^s v^{H-\frac{1}{2}}(v-t)^{H-\frac{3}{2}}(s-t)^{\frac{1}{2}-H}dv\right)^2\right)$$

$$\leq\frac{2c_H^2\left(\left(B(H-\frac{1}{2},\frac{3}{2}-H)\right)^2+\frac{1}{(H-\frac{1}{2})^2}\right)t^{1-2H}}{(\frac{1}{2}-H)^2},$$

where $B(\cdot,\cdot)$ is beta function. By (10) and (12),

$$\int_t^1\left(\int_t^s\frac{(s^{\frac{1}{2}-H}-u^{\frac{1}{2}-H})K(u,t)}{(s-u)^{\frac{1}{2}+H}}du s^{H-\frac{1}{2}}\right)^2 ds\leq\int_t^1 C_1^2 C_2^2 s^{2-4H}t^{1-2H}s^{2H-1}ds\leq\frac{C_1^2 C_2^2}{2-2H}t^{1-2H}. \quad (42)$$

Then, combining (39), (40), (41) and (42) yields

$$\left(\int_t^1 K(s,t)P_s ds\right)^2$$

$$\leq\frac{C_1^2 t^{1-2H}}{\Gamma^2(\frac{3}{2}-H)(2-2H)}\int_0^1 |(K^{-1}DF)_s|^2 ds+\frac{4C_1^2 t^{1-2H}}{\Gamma^2(\frac{3}{2}-H)(2-2H)}\int_0^1 |(K^{-1}DF)_s|^2 ds$$

$$+\frac{2c_H^2\left(\left(B(H-\frac{1}{2},\frac{3}{2}-H)\right)^2+\frac{1}{(H-\frac{1}{2})^2}\right)t^{1-2H}}{\Gamma^2(\frac{3}{2}-H)}\int_0^1 |(K^{-1}DF)_s|^2 ds$$

$$+\frac{C_1^2 C_2^2(H-\frac{1}{2})^2 t^{1-2H}}{\Gamma^2(\frac{3}{2}-H)(2-2H)}\int_0^1 |(K^{-1}DF)_s|^2 ds.$$

It follows that

$$\left(\int_t^1 K(s,t)P_s ds\right)^2 dt\leq\widehat{C}t^{1-2H}\int_0^1 |(K^{-1}DF)_s|^2 ds, \quad (43)$$

where

$$\widehat{C}=\frac{C_1^2}{\Gamma^2(\frac{3}{2}-H)(2-2H)}+\frac{4C_1^2}{\Gamma^2(\frac{3}{2}-H)(2-2H)}$$

$$+\frac{2c_H^2\left(\left(B(H-\frac{1}{2},\frac{3}{2}-H)\right)^2+\frac{1}{(H-\frac{1}{2})^2}\right)}{\Gamma^2(\frac{3}{2}-H)}+\frac{C_1^2 C_2^2(H-\frac{1}{2})^2}{\Gamma^2(\frac{3}{2}-H)(2-2H)}.$$

Therefore, by (30), (36) and (43),

$$
\begin{aligned}
\mathbb{E}_\nu \left[ \int_0^1 \frac{|\eta_t|^2}{G_t} dt \right] \le & 8\mathbb{E}_\nu \left[ \int_0^1 |(K^{-1}DF)_t|^2 dt + 4\alpha^4 e^{2\alpha} C_1^2 (1+\alpha^2) \int_0^1 t^{1-2H} dt C \int_0^1 |(K^{-1}DF)_s|^2 ds \right. \\
& \left. + 2\alpha^2 \widehat{C} \int_0^1 t^{1-2H} dt \int_0^1 |(K^{-1}DF)_s|^2 ds \right] \\
= & 8 \left( 1 + \frac{4\alpha^4 e^{2\alpha} C_1^2 (1+\alpha^2) C}{2-2H} + \frac{2\alpha^2 \widehat{C}}{2-2H} \right) \mathbb{E}_\nu \left[ \int_0^1 |(K^{-1}DF)_s|^2 ds \right].
\end{aligned}
$$

Hence, by (27), we obtain the Logarithmic-Sobolev inequality

$$
\mathbb{E}_\nu [F^2 \ln F^2] \le 4 \left( 1 + \frac{4\alpha^4 e^{2\alpha} C_1^2 (1+\alpha^2) C}{2-2H} + \frac{2\alpha^2 \widehat{C}}{2-2H} \right) \mathbb{E}_\nu \left[ \int_0^1 |(K^{-1}DF)_s|^2 ds \right] + \mathbb{E}_\nu [F^2] \ln \mathbb{E}_\nu [F^2].
$$

Now we complete the proof. $\qquad\square$

**Acknowledgments.** Xiaoxia Sun is supported by the Research Foundation of Dongbei University of Finance and Economics under grants 2015104, the Postdoctoral Science Foundation of Dongbei University of Finance and Economics under grants BSH201519. Feng Guo is supported by the National Natural Science Foundation of China under grants 11401074.

## References

[1] S. Aida. Differential calculus on path and loop spaces II: Irreducibility of dirichlet forms on loop spaces. *Bulletin des Sciences Mathématiques*, 122(8):635–666, 1998.

[2] S. Aida. Logarithmic derivatives of heat kernels and Logarithmic Sobolev inequalities with unbounded diffusion coefficients on loop spaces. *Journal of Functional Analysis*, 174:430–477, 2000.

[3] S. Aida and K. D. Elworthy. Differential calculus on path and loop spaces, I. Logarithmic Sobolev inequalities on path spaces. *Comptes Rendus de l'Académie des Sciences (Paris) Serié I*, 321:97–102, 1995.

[4] J. M. Bismut. *Large Deviations and the Malliavin Calculus*. Birkhäuser, 1984.

[5] M. Capitaine, E. P. Hsu, and M. Ledoux. Martingale representation and a simple proof of Logarithmic Sobolev inequalities on path spaces. *Electronic Communications in Probability*, 2:71–81, 1997.

[6] L. Decreusefond and A. S. Üstünel. Stochastic analysis of the fractional Brownian motion. *Potential Analysis*, 10(2):177–214, 1999.

[7] B. K. Driver. A Cameron-Martin type quasi-invariance theorem for Brownian motion on a compact Riemannian manifold. *Journal of Functional Analysis*, 110(2):272–376, 1992.

[8] B. K. Driver. A Cameron-Martin type quasi-invariance theorem for pinned Brownian motion on a compact Riemannian manifold. *Transactions of the American Mathematical Society*, 342(1):375–395, 1994.

[9] T. E. Duncan, Y. Z. Hu, and B. Pasik-Duncan. Stochastic calculus for fractional Brownian motion I. Theory. *SIAM Journal on Control and Optimization*, 38(2):582–612, 2000.

[10] O. Enchev and D. W. Stroock. Pinned Brownian motion and its perturbations. *Advances in Mathematics*, 119(2):127–154, 1996.

[11] S. Z. Fang and P. Malliavin. Stochastic analysis on the path space of a Riemannian manifold: I. Markovian stochastic calculus. *Journal of Functional Analysis*, 118(1):249–274, 1993.

[12] F. Z. Gong and Z. M. Ma. The Log-Sobolev inequality on loop space over a compact Riemannian manifold. *Journal of Functional Analysis*, 157:599–623, 1998.

[13] L. Gross. Logarithmic Sobolev inequalities. *American Journal of Mathematics*, 97:1061–1083, 1975.

[14] L. Gross. Logarithmic Sobolev inequalities on loop groups. *Journal of Functional Analysis*, 102:268–313, 1991.

[15] E. P. Hsu. Quasi-invariance of the Wiener measure on the path space over a compact Riemannian manifold. *Journal of Functional Analysis*, 134(2):417–450, 1995.

[16] E. P. Hsu. Integration by parts in loop spaces. *Mathematische Annalen*, 309(2):331–339, 1997.

[17] E. P. Hsu. Logarithmic Sobolev inequalities on path spaces over Riemannian manifolds. *Communications in Mathematical Physics*, 189:9–16, 1997.

[18] E. P. Hsu. *Stochastic Analysis on Manifolds*. American Mathematical Society, 2002.

[19] D. Nualart and Y. Ouknine. Regularization of differential equations by fractional noise. *Stochastic Processes and their Applications*, 102(1):103–116, 2002.

[20] S. G. Samko, A. A. Kilbas, and O. I. Marichev. *Fractional Integrals and Derivatives: Theory and Applications*. Gordon Breach Science, 1993.

[21] X. X. Sun and F. Guo. On integration by parts formula and characterization of fractional Ornstein-Uhlenbeck process. *Statistics and Probability Letters*, 107(5):170–177, 2015.

14